# On the Fontaine-Mazur Conjecture for CM-Fields

by Kay Wingberg at Heidelberg

In [3] Fontaine and Mazur conjecture (as a consequence of a general principle) that a number field $k$ has no infinite unramified Galois extension such that its Galois group is a $p$-adic analytic pro-$p$-group. A counter-example to this conjecture would produce an unramified Galois representation with infinite image, that could not "come from geometry". Some evidence for this conjecture is shown in [1] and [4].

Since every $p$-adic analytic pro-$p$-group contains an open powerful resp. uniform subgroup one is led to the question whether a given number field possesses an infinite unramified Galois $p$-extension with powerful resp. uniform Galois group. With regard to this problem, we would like to mention the main result of Boston, [1] theorem 1:

> Let $p$ be a prime number and let $k|k_0$ be a finite cyclic Galois extension of degree prime to $p$ such that $p$ does not divide the class number of $k_0$. Then, if the Galois group $G(M|k)$ of an unramified Galois $p$-extension $M$ of $k$ is powerful, it is finite.

In this paper we will prove a statement which is in some sense weaker as the above and in another sense stronger (and in view of the general conjecture very weak):

> Let $p$ be odd and let $k$ be a CM-field with maximal totally real subfield $k^+$ containing the group $\mu_p$ of $p$-th roots of unity. Let $M = L(p)$ be the maximal unramified $p$-extension of $k$. Assume that the $p$-rank of the ideal class group $Cl(k^+)$ of $k^+$ is not equal to 1. Then, if the Galois group $G(L(p)|k)$ is powerful, it is finite.

If the $p$-rank of $Cl(k^+)$ is equal to 1, we have two weaker results. First, replacing the word powerful by uniform and assuming that the first step in the $p$-cyclotomic tower of $k$ is not unramified, then the statement above holds without any condition on $Cl(k^+)$. Secondly, we consider the conjecture in the $p$-cyclotomic tower of the number field $k$. Denote the $n$-th layer of the cyclotomic $\mathbb{Z}_p$-extension $k_\infty$ of $k$ by $k_n$ and let $G(L_n(p)|k_n)$ be the Galois group of the maximal unramified $p$-extension $L_n(p)$ of $k_n$. Then the following statement holds.

> Let $p \neq 2$ and let $k$ be a CM-field containing $\mu_p$. Assume that the Iwasawa $\mu$-invariant of $k_\infty|k$ is zero. Then there exists a number $n_0$ such that for all $n \geq n_0$ the following holds: If the Galois group $G(L_n(p)|k_n)$ is powerful, then it is finite.

Let $S$ be a set of primes of $k$ containing the set $S_\infty$ of archimedean primes and assume that no prime of $S$ split in the extension $k|k^+$. Then all the results above hold, if we replace the field $L(p)$ by the maximal unramified $p$-extension $L_S(p)$ which is completely decomposed at all primes in $S$ and the ideal class group $Cl(k^+)$ by the $S$-ideal class group $Cl_S(k^+)$ of $k^+$.

Of course, our main interest is the conjecture for general $p$-adic analytic groups. We will prove the following result.

> Let $p \neq 2$ and let $k$ be a CM-field containing $\mu_p$ with maximal totally real subfield $k^+$ and assume that $\mu_p \not\subseteq k^+_{\mathfrak{p}}$ for all primes $\mathfrak{p}$ of $k^+$ above $p$. Then, if $G(L_k(p)|k)$ is $p$-adic analytic, $G(L_{k^+}(p)|k^+)$ is finite.

Unfortunately, we do not have Boston's result for general analytic pro-$p$-groups. Otherwise, in the situation above it would follow that $G(L_k(p)|k)$ is not an infinite $p$-adic analytic group.

# 1 A duality theorem

We use the following notation:

| | |
|---|---|
| $p$ | is a prime number, |
| $k$ | is a number field, |
| $S_\infty$ | is the set of archimedean primes of $k$, |
| $S$ | is a set of primes of $k$ containing $S_\infty$, |
| $E_S(k)$ | is the group of $S$-units of $k$, |
| $Cl_S(k)$ | is the $S$-ideal class group of $k$, |
| $L_S$ | is the maximal unramified extension of $k$ which is completely decomposed at $S$, |
| $L_S(p)$ | is the maximal $p$-extension of $k$ inside $L_S$, |
| $L$ | is the maximal unramified extension of $k$, |
| $L(p)$ | is the maximal $p$-extension of $k$ inside $L$. |

We write $E(k)$ for the group $E_{S_\infty}(k)$ of units of $k$ and $Cl(k)$ for the ideal class group $Cl_{S_\infty}(k)$ of $k$. Obviously,

$$\begin{aligned} L &= L_{S_\infty}, & \text{if } k \text{ is totally imaginary,} \\ L(p) &= L_{S_\infty}(p), & \text{if } p \neq 2 \text{ or } k \text{ totally imaginary.} \end{aligned}$$

If $K$ is an infinite algebraic extension of $\mathbb{Q}$, then $E_S(K) = \varinjlim_k E_S(k)$ where $k$ runs through the finite subextensions of $K$.

For a profinite group $G$, a discrete $G$-module $M$ and any integer $i$ the $i$-th Tate cohomology is defined by

$$\hat{H}^i(G,M) = H^i(G,M) \text{ for } i \geq 1 \text{ and } \hat{H}^i(G,M) = \varprojlim_{U, def} \hat{H}^i(G/U, M^U) \text{ for } i \leq 0,$$

where $U$ runs through all open normal subgroups of $G$ and the transition maps are given by the deflation, see [7].

**Theorem 1.1** *Let $S$ be a set of primes of $k$ containing $S_\infty$. Then the following holds:*

(i) *There are canonical isomorphisms*

$$\hat{H}^i(G(L_S|k), E_S(L_S)) \cong \hat{H}^{2-i}(G(L_S|k), \mathbb{Q}/\mathbb{Z})^\vee$$

*for all $i \in \mathbb{Z}$. Here $^\vee$ denotes the Pontryagin dual.*

(ii) *There are canonical isomorphisms*

$$\hat{H}^i(G(L_S(p)|k), E_S(L_S(p))) \cong \hat{H}^{2-i}(G(L_S(p)|k), \mathbb{Q}_p/\mathbb{Z}_p)^\vee$$

*for all $i \in \mathbb{Z}$.*

**Proof:** Let $C_S(L_S)$ be the $S$-idele class group of $L_S$. The subgroup $C_S^0(L_S)$ of $C_S(L_S)$ given by the ideles of norm 1 is a level-compact class formation for $G(L_S|k)$ with divisible group of universal norms. From the duality theorem of Nakayama-Tate we obtain the isomorphisms

$$\hat{H}^i(G(L_S|k), C_S(L_S)) \cong \hat{H}^{2-i}(G(L_S|k), \mathbb{Z})^\vee, \quad i \in \mathbb{Z},$$

since $\hat{H}^i(G(L_S|k), C_S(L_S)) \cong \hat{H}^i(G(L_S|k), C_S^0(L_S))$, see [7] proposition 4.

Let $K|k$ be a finite Galois extension inside $L_S$. From the exact sequence

$$0 \longrightarrow E_S(K) \longrightarrow J_S(K) \longrightarrow C_S(K) \longrightarrow Cl_S(K) \longrightarrow 0,$$

where $J_S(K)$ denotes the group of $S$-ideles of $K$, which is a cohomological trivial $G(K|k)$-module ($K|k$ is completely decomposed at $S$), we obtain isomorphisms

$$\hat{H}^{i+1}(G(K|k), E_S(K)) \cong \hat{H}^i(G(K|k), D(K)),$$

where $D(K)$ denotes the kernel of the surjection $C_S(K) \twoheadrightarrow Cl_S(K)$, and a long exact sequences

$$\longrightarrow \hat{H}^i(G(K|k), D(K)) \longrightarrow \hat{H}^i(G(K|k), C_S(K)) \longrightarrow \hat{H}^i(G(K|k), Cl_S(K)) \longrightarrow .$$

If $K'$ is the maximal abelian extension of $K$ in $L_S$, then $G(L_S|K')$ is an open subgroup of $G(L_S|K)$ by the finiteness of the class number of $K$. The commutative diagram

$$\begin{array}{ccc} Cl_S(K') & \xrightarrow{norm} & Cl_S(K) \\ {\scriptstyle rec}\downarrow\wr & & {\scriptstyle rec}\downarrow\wr \\ G(L_S|K')^{ab} & \xrightarrow{can} & G(L_S|K)^{ab} \end{array}$$

shows, since $can$ is the zero map, that

$$Cl_S(K') \xrightarrow{norm} Cl_S(K)$$

is trivial. It follows that

$$\varprojlim_K \hat{H}^i(G(K|k), Cl_S(K)) = 0 \quad \text{for } i \leq 0.$$

Since all groups in the exact sequence above are finite, we can pass to the projective limit and we obtain isomorphisms

$$\varprojlim_K \hat{H}^i(G(K|k), D(K)) \cong \hat{H}^i(G(L_S|k), C_S(L_S)) \quad \text{for } i \leq 0,$$

and therefore isomorphisms

$$\hat{H}^{i+1}(G(L_S|k), E_S(L_S)) \cong \hat{H}^i(G(L_S|k), C_S(L_S)) \quad \text{for } i \leq -1.$$

The last assertion also holds for $i = 0$: from the commutative diagram

$$\begin{array}{ccc} \hat{H}^0(G(K'|k), D(K')) & \xrightarrow{\delta}_{\sim} & H^1(G(K'|k), E_S(K')) \\ \downarrow{\scriptstyle def} & & \downarrow \\ \hat{H}^0(G(K|k), D(K)) & \xrightarrow{\delta}_{\sim} & H^1(G(K|k), E_S(K)), \end{array}$$

where $k \subseteq K \subseteq K'$ are finite Galois extensions inside $L_S$, it follows that the limit $\varprojlim_K H^1(G(K|k), E_S(K))$ exists. Since

$$H^1(G(K|k), E_S(K)) \subseteq H^1(G(L_S|k), E_S(L_S)) \cong Cl_S(k)$$

and

$$\begin{aligned} \varprojlim_K \hat{H}^0(G(K|k), D(K)) &\cong \hat{H}^0(G(L_S|k), C_S(L_S)) \cong H^2(G(L_S|k), \mathbb{Z})^\vee \\ &\cong H^1(G(L_S|k), \mathbb{Q}/\mathbb{Z})^\vee = G(L_S|k)^{ab} \cong Cl_S(k), \end{aligned}$$

the projective limit $\varprojlim_K H^1(G(K|k), E_S(K))$ becomes stationary and is equal to $H^1(G(L_S|k), E_S(L_S))$.

For $i \geq 1$ the exact sequence

$$0 \longrightarrow E_S(L_S) \longrightarrow J_S(L_S) \longrightarrow C_S(L_S) \longrightarrow 0$$

induces isomorphisms

$$H^i(G(L_S|k), C_S(L_S)) \cong H^{i+1}(G(L_S|k), E_S(L_S)).$$

Putting all together, we obtain canonical isomorphisms

$$\hat{H}^{i+1}(G(L_S|k), E_S(L_S)) \cong \hat{H}^{2-i}(G(L_S|k), \mathbb{Z})^\vee \cong \hat{H}^{1-i}(G(L_S|k), \mathbb{Q}/\mathbb{Z})^\vee$$

for all $i \in \mathbb{Z}$. The proof for the field $L_S(p)$ is analogously. $\square$

Let $k$ be a number field of CM-type with maximal totally real subfield $k^+$ and let $\Delta = G(k|k^+) = \langle \sigma \rangle \cong \mathbb{Z}/2\mathbb{Z}$. If $p \neq 2$, we put as usual

$$M^\pm = (1 \pm \sigma)M$$

for a $\mathbb{Z}_p[\Delta]$-module $M$. For a $\mathbb{Z}_p$-module $N$ let $_pN = \{x \in N \mid px = 0\}$.

**Corollary 1.2** *Let $p$ be an odd prime number and let $k$ be a CM-field. Let $S$ be a set of primes of $k$ containing $S_\infty$ and assume that no prime of $S$ split in the extension $k|k^+$. Then*

$$\dim_{\mathbb{F}_p} {}_pH^2(G(L_S(p)|k), \mathbb{Q}_p/\mathbb{Z}_p)^- \leq \delta,$$

*where $\delta$ is equal to 1 if $k$ contains the group $\mu_p$ of $p$-th roots of unity and otherwise equal to 0.*

**Proof:** By proposition 1.1, there is a $\Delta$-invariant surjection

$$E_S(k) \twoheadrightarrow \hat{H}^0(G(L_S(p)|k), E_S(L_S(p))) \cong H^2(G(L_S(p)|k), \mathbb{Q}_p/\mathbb{Z}_p)^\vee$$

and so a surjection

$$(E_S(k)/p)^- \twoheadrightarrow ({}_pH^2(G(L_S(p)|k), \mathbb{Q}_p/\mathbb{Z}_p)^-)^\vee.$$

Since no prime of $S$ splits in the extension $k|k^+$, we have $(E_S(k)/p)^- \cong \mu_p(k)$ which gives us the desired result. $\square$

# 2 Powerful pro-$p$-groups with involution

Let $p$ be a prime number. For a pro-$p$-group $G$ the descending $p$-central series is defined by
$$G_1 = G, \qquad G_{i+1} = (G_i)^p[G_i, G] \quad \text{for } i \geq 1.$$
If a group $\Delta \cong \mathbb{Z}/2\mathbb{Z}$ acts on $G$ and $p$ is odd, then we define
$$d(G)^\pm = \dim_{\mathbb{F}_p}(G/G_2)^\pm = \dim_{\mathbb{F}_p} H^1(G, \mathbb{Z}/p\mathbb{Z})^\pm.$$

The following proposition also follows from Boston result (resp. its proof), but in our situation, where only an involution acts on $G$, we will give a simple proof.

**Proposition 2.1** *Let $p \neq 2$ and let $G$ be a finitely generated powerful pro-$p$-group with an action by the group $\Delta \cong \mathbb{Z}/2\mathbb{Z}$. Then the following holds:*

*If $d(G)^+ = 0$, then $G$ is abelian.*

*In particular, if $d(G)^+ = 0$ and $G^{ab}$ is finite, then $G$ is finite.*

**Proof:** Since $G$ is powerful, we have
$$[G, G]/H \subseteq G^p H/H \quad \text{where } H = ([G, G])^p[G, G, G].$$
From $G/G_2 = (G/G_2)^-$ it follows that
$$[G, G]/H = ([G, G]/H)^+ \quad \text{and} \quad G^p H/H = (G^p H/H)^-,$$
since $G/[G, G] = (G/[G, G])^-$ and $G^p = \{x^p \mid x \in G\}$, [2] theorem 3.6(iii), and so
$$(x^p)^\sigma \equiv x^{-p} \mod H \quad \text{for } 1 \neq \sigma \in \Delta \text{ and } x \in G.$$
We obtain
$$[G, G] \subseteq ([G, G])^p[G, G, G].$$
This implies $[G, G] = 1$. $\square$

**Proposition 2.2** *Let $p \neq 2$ and let $G$ be a finitely generated powerful pro-$p$-group with an action by the group $\Delta \cong \mathbb{Z}/2\mathbb{Z}$. Assume that $G^{ab}$ is finite. Then the following inequalities hold:*

(i) $\quad d(G)^+ \cdot d(G)^- \leq d(G)^- + \dim_{\mathbb{F}_p} {}_p H^2(G, \mathbb{Q}_p/\mathbb{Z}_p)^-,$

(ii) $\quad \binom{d(G)^+}{2} + \binom{d(G)^-}{2} \leq d(G)^+ + \dim_{\mathbb{F}_p} {}_p H^2(G, \mathbb{Q}_p/\mathbb{Z}_p)^+.$

**Proof:** Let $d^\pm = d(G)^\pm$. From the exact sequences

$$0 \longrightarrow H^1(G/G_2, \mathbb{Z}/p\mathbb{Z}) \xrightarrow{\sim} H^1(G, \mathbb{Z}/p\mathbb{Z}) \longrightarrow H^1(G_2, \mathbb{Z}/p\mathbb{Z})^G$$
$$\longrightarrow H^2(G/G_2, \mathbb{Z}/p\mathbb{Z}) \longrightarrow H^2(G, \mathbb{Z}/p\mathbb{Z})$$

and

$$0 \longrightarrow (_pG^{ab})^\vee \longrightarrow H^2(G, \mathbb{Z}/p\mathbb{Z}) \longrightarrow {}_pH^2(G, \mathbb{Q}_p/\mathbb{Z}_p) \longrightarrow 0$$

we obtain the inequalities

$$\dim_{\mathbb{F}_p} H^2(G/G_2, \mathbb{Z}/p\mathbb{Z})^\pm \le \dim_{\mathbb{F}_p}(G_2/G_3)^\pm + d^\pm + \dim_{\mathbb{F}_p} {}_pH^2(G, \mathbb{Q}_p/\mathbb{Z}_p)^\pm.$$

Here we used $\dim_{\mathbb{F}_p}(_pG^{ab})^\pm = d^\pm$ which holds by the finiteness of $G^{ab}$. Since $G$ is powerful, the $\Delta$-invariant homomorphism

$$G/G_2 \xrightarrow{p} G_2/G_3$$

is surjective, see [2] theorem 3.6, and we obtain

$$\dim_{\mathbb{F}_p} H^2(G/G_2, \mathbb{Z}/p\mathbb{Z})^\pm \le 2d^\pm + \dim_{\mathbb{F}_p} {}_pH^2(G, \mathbb{Q}_p/\mathbb{Z}_p)^\pm.$$

Let

$$G/G_2 \cong A_1 \oplus \cdots \oplus A_{d^+} \oplus B_1 \oplus \cdots \oplus B_{d^-}$$

be a $\Delta$-invariant decomposition into cyclic groups of order $p$ such that $A_i = A_i^+$ and $B_j = B_j^-$. For $H^2(G/G_2, \mathbb{Z}/p\mathbb{Z})$ we obtain the $\Delta$-invariant Künneth decomposition:

$$\begin{aligned}
H^2(G/G_2, \mathbb{Z}/p\mathbb{Z}) &\cong \bigoplus_{i=1}^{d^+} H^2(A_i, \mathbb{Z}/p\mathbb{Z}) \\
&\oplus \bigoplus_{i<j} H^1(A_i, \mathbb{Z}/p\mathbb{Z}) \otimes H^1(A_j, \mathbb{Z}/p\mathbb{Z}) \\
&\oplus \bigoplus_{i<j} H^1(B_i, \mathbb{Z}/p\mathbb{Z}) \otimes H^1(B_j, \mathbb{Z}/p\mathbb{Z}) \\
&\oplus \bigoplus_{i=1}^{d^-} H^2(B_i, \mathbb{Z}/p\mathbb{Z}) \\
&\oplus \bigoplus_{i,j} H^1(A_i, \mathbb{Z}/p\mathbb{Z}) \otimes H^1(B_j, \mathbb{Z}/p\mathbb{Z}).
\end{aligned}$$

Counting dimensions yields

$$\begin{aligned}
\dim_{\mathbb{F}_p} H^2(G/G_2, \mathbb{Z}/p\mathbb{Z})^+ &= d^+ + \binom{d^+}{2} + \binom{d^-}{2}, \\
\dim_{\mathbb{F}_p} H^2(G/G_2, \mathbb{Z}/p\mathbb{Z})^- &= d^- + d^+ d^-.
\end{aligned}$$

This proves the proposition. $\square$

Now we analyze the case where $G$ is a powerful pro-$p$-group which is a Poincaré group of dimension 3.

**Proposition 2.3** *Let $p$ be odd and let $P$ be a finitely generated powerful pro-$p$-group with an action of $\Delta \cong \mathbb{Z}/2\mathbb{Z}$.*

(i) *If $P$ is uniform, then*
$$\dim_{\mathbb{F}_p} H^2(P, \mathbb{Z}/p\mathbb{Z})^+ = \binom{d(P)^+}{2} + \binom{d(P)^-}{2},$$
$$\dim_{\mathbb{F}_p} H^2(P, \mathbb{Z}/p\mathbb{Z})^- = d(P)^+ \cdot d(P)^-.$$

(ii) *If $P$ is uniform such that $P^{ab}$ is finite and $d(P)^+ = 1$, then*
$$\dim_{\mathbb{F}_p} {}_p H^2(P, \mathbb{Q}_p/\mathbb{Z}_p)^- = 0.$$

(iii) *If $P$ is a Poincaré group of dimension $3$ such that $P^{ab}$ is finite, then*
$$d(P)^+ = 1 \quad and \quad d(P)^- = 2 \quad\quad or$$
$$d(P)^+ = 3 \quad and \quad d(P)^- = 0.$$

**Proof:** Let $P$ be uniform. By [2] definition 4.1 and theorem 4.26, we have
$$\dim_{\mathbb{F}_p}(H^1(P_2, \mathbb{Z}/p\mathbb{Z})^P)^\pm = d(P)^\pm \quad \text{and} \quad \dim_{\mathbb{F}_p} H^2(P, \mathbb{Z}/p\mathbb{Z}) = \binom{d(P)}{2}.$$

Counting dimensions shows that
$$\dim_{\mathbb{F}_p} H^2(P/P_2, \mathbb{Z}/p\mathbb{Z}) = \dim_{\mathbb{F}_p} H^1(P_2, \mathbb{Z}/p\mathbb{Z})^P + \dim_{\mathbb{F}_p} H^2(P, \mathbb{Z}/p\mathbb{Z}),$$

and so the sequence
$$0 \longrightarrow H^1(P_2, \mathbb{Z}/p\mathbb{Z})^P \longrightarrow H^2(P/P_2, \mathbb{Z}/p\mathbb{Z}) \longrightarrow H^2(P, \mathbb{Z}/p\mathbb{Z}) \longrightarrow 0$$

is exact. Therefore
$$\dim_{\mathbb{F}_p} H^2(P, \mathbb{Z}/p\mathbb{Z})^\pm = \dim_{\mathbb{F}_p} H^2(P/P_2, \mathbb{Z}/p\mathbb{Z})^\pm - \dim_{\mathbb{F}_p}(H^1(P_2, \mathbb{Z}/p\mathbb{Z})^P)^\pm,$$

which proves (i).

If $P^{ab}$ is finite and $d(P)^+ = 1$, then by (i)
$$\dim_{\mathbb{F}_p} {}_p H^2(P, \mathbb{Q}_p/\mathbb{Z}_p)^- = \dim_{\mathbb{F}_p} H^2(P, \mathbb{Z}/p\mathbb{Z})^- - \dim_{\mathbb{F}_p}({}_p P^{ab})^-$$
$$= d(P)^+ \cdot d(P)^- - d(P)^- = 0.$$

Now let $P$ be a powerful Poincaré group of dimension 3; in particular, $P$ is torsionfree and therefore $P$ is uniform, see [2] theorem 4.8. Since
$$\dim_{\mathbb{F}_p} H^1(P, \mathbb{Z}/p\mathbb{Z}) = \dim_{\mathbb{F}_p} H^2(P, \mathbb{Z}/p\mathbb{Z})$$

and since $P^{ab}$ is finite, the exact sequence

$$0 \longrightarrow ({}_pP^{ab})^\vee \longrightarrow H^2(P, \mathbb{Z}/p\mathbb{Z}) \longrightarrow {}_pH^2(P, \mathbb{Q}_p/\mathbb{Z}_p) \longrightarrow 0$$

shows that

$$({}_pP^{ab})^\vee \xrightarrow{\sim} H^2(P, \mathbb{Z}/p\mathbb{Z}).$$

It follows that

$$\dim_{\mathbb{F}_p} H^2(P, \mathbb{Z}/p\mathbb{Z})^\pm = d(P)^\pm,$$

and so by (i)

$$d(P)^+ \cdot d(P)^- = d(P)^-.$$

This proves (iii). □

## 3 On the Fontaine-Mazur Conjecture

We keep the notation of sections 1 and 2. Let

$$d_k^\pm = \dim_{\mathbb{F}_p}(Cl(k)/p)^\pm = d(G(L(p)|k))^\pm.$$

**Theorem 3.1** *Let $p$ be an odd prime number and let $k$ be a CM-field such that*
(i) $d_k^- \neq 0$, if $\mu_p \not\subseteq k$,
(ii) $d_k^+ \neq 1$.
*Then, if the Galois group $G(L(p)|k)$ of the maximal unramified $p$-extension $L(p)$ of $k$ is powerful, it is finite.*

**Proof:** If $d_k^+ = 0$, then the theorem follows from proposition 2.1. Therefore we assume that $d_k^+ \geq 2$ (assumption (ii)). From assumption (i) and Leopoldt's Spiegelungssatz, see [8] theorem 10.11, it follows that $d_k^- \geq 1$. From proposition 2.2 and corollary 1.2 we obtain the inequality

$$d_k^+ d_k^- \leq d_k^- + \delta.$$

It follows that $d_k^+ = 2$, $d_k^- = 1$.

Let $P = G(L(p)|k)_i$, $i$ large enough. Then $P$ is uniform, [2] theorem 4.2, and $d(P) \leq 3$, [2] theorem 3.8. Furthermore, if $P$ is non-trivial, then $P$ is a Poincaré group of dimension $\dim(P) = d(P) \leq 3$, see [5] chap.V theorem (2.2.8) and (2.5.8). But Poincaré groups of dimension $\dim(P) \leq 2$ have the group $\mathbb{Z}_p$ as homomorphic image, and so we can assume that $\dim(P) = d(P) = 3$. Since $G(L(p)|k)$ is powerful, we have a surjection

$$G(L(p)|k)/G(L(p)|k)_2 \twoheadrightarrow G(L(p)|k)_i/G(L(p)|k)_{i+1}.$$

Furthermore, by [2] theorem 3.6(ii), $G(L(p)|k)_{i+1} = (G(L(p)|k)_i)_2 = P_2$, and so $G(L(p)|k)_i/G(L(p)|k)_{i+1} = P/P_2$. Therefore $d(P)^+ = 2$ and $d(P)^- = 1$. Now the result is a consequence of proposition 2.3(iii). □

If $\mu_p \subseteq k$, then $d_k^+ = 1$ is the only remaining case. Here we only get a weaker result. Let $k_\infty$ be the cyclotomic $\mathbb{Z}_p$-extension of $k$ and denote by $k_n$ the $n$-th layer of $k_\infty|k$.

**Theorem 3.2** *Let $p \neq 2$ and let $k$ be a CM-field containing $\mu_p$. Assume that $k_1|k$ is not unramified if $d_k^+ = 1$. Then the Galois group $G(L(p)|k)$ of the maximal unramified $p$-extension $L(p)$ of $k$ is not uniform.*

**Proof:** Suppose that $G = G(L(p)|k)$ is uniform. Using theorem 3.1, we may assume that $d(G)^+ = 1$, and so, by proposition 2.3(ii),
$$\dim_{\mathbb{F}_p} {}_p H^2(G, \mathbb{Q}_p/\mathbb{Z}_p)^- = 0.$$

On the other hand, by theorem 1.1, we have a surjection
$$H^2(G, \mathbb{Q}_p/\mathbb{Z}_p)^\vee \cong \hat{H}^0(G, E_{L(p)}) \twoheadrightarrow \hat{H}^0(G(K|k), E_K)$$
where $K|k$ is a finite unramified Galois $p$-extension of CM-fields (recall that $d(G)^+ \neq 0$), and so a surjection
$$(H^2(G, \mathbb{Q}_p/\mathbb{Z}_p)^-)^\vee \twoheadrightarrow \hat{H}^0(G(K|k), E_K)^-.$$

Since $K$ is of CM-type, it follows that
$$\hat{H}^0(G(K|k), E_K)^- \cong \hat{H}^0(G(K|k), \mu(K)(p)).$$

By our assumption, $K$ is disjoint to $k_\infty$, i.e. $\mu(K)(p) = \mu(k)(p)$, and so
$$\dim_{\mathbb{F}_p} \hat{H}^0(G(K|k), \mu(K)(p)) = 1.$$

It follows that
$$\dim_{\mathbb{F}_p} {}_p H^2(G, \mathbb{Q}_p/\mathbb{Z}_p)^- = 1.$$

This contradiction proves the theorem. □

Now we consider the Galois groups $G(L_n(p)|k_n)$ of the maximal unramified $p$-extension $L_n(p)$ of $k_n$ in the $p$-cyclotomic tower of $k$.

**Theorem 3.3** *Let $p \neq 2$ and let $k$ be a CM-field containing $\mu_p$. Assume that the Iwasawa $\mu$-invariant of the cyclotomic $\mathbb{Z}_p$-extension $k_\infty|k$ is zero.*

*Then there exists a number $n_0$ such that for all $n \geq n_0$ the following holds: If the Galois group $G(L_n(p)|k_n)$ is powerful, then it is finite.*

**Proof:** Let
$$1 \longrightarrow G_\infty \longrightarrow G(L_\infty(p)|k) \longrightarrow \Gamma \longrightarrow 1$$
where $G_\infty = G(L_\infty(p)|k_\infty)$ is the Galois group of the maximal unramified $p$-extension $L_\infty(p)$ of $k_\infty$ and $\Gamma = G(k_\infty|k) = \langle \gamma \rangle$. Let $\Gamma_n = \langle \gamma^{p^n} \rangle$, $n \geq 0$, be the open subgroups of $\Gamma$ of index $p^n$. By our assumption on the Iwasawa $\mu$-invariant $G_\infty$ is a finitely generated pro-$p$-group.

Let $n_1$ be large enough such that all primes of $k_{n_1}$ above $p$ are totally ramified in $k_\infty|k_{n_1}$ and let $\langle \gamma_j \rangle \subseteq G(k_\infty|k_{n_1})$, $j = 1, \ldots, s$, be the inertia groups of some extensions of the finitely many primes $\mathfrak{p}_1, \ldots \mathfrak{p}_s$ of $k_{n_1}$ above $p$.

For $n \geq n_1$ let
$$M_n = (\gamma_j^{p^{n-n_1}}, j = 1, \ldots, s) \subseteq G(L_\infty(p)|k_n)$$
be the normal subgroup generated by all conjugates of the elements $\gamma_j^{p^{n-n_1}}$ and
$$N_n := M_n \cap G_\infty = (\gamma_i^{p^{n-n_1}} \gamma_j^{-p^{n-n_1}}, [\gamma_j^{p^{n-n_1}}, g], i,j = 1, \ldots, s, \, g \in G_\infty).$$

Then the commutative exact diagram

$$\begin{array}{ccccccccc}
1 & \longrightarrow & N_n & \longrightarrow & M_n & \longrightarrow & \Gamma_n & \longrightarrow & 1 \\
& & \downarrow & & \downarrow & & \| & & \\
1 & \longrightarrow & G_\infty & \longrightarrow & G(L_\infty(p)|k_n) & \longrightarrow & \Gamma_n & \longrightarrow & 1
\end{array}$$

shows that
$$G_\infty/N_n \cong G(L_n(p)|k_n)$$
and we have canonical surjections
$$G_\infty \twoheadrightarrow G(L_m(p)|k_m) \twoheadrightarrow G(L_n(p)|k_n)$$
for $m \geq n \geq n_1$.

Let $n_0 \geq n_1$ be large enough such that
$$G_\infty/(G_\infty)_3 \xrightarrow{\sim} G(L_n(p)|k_n)/(G(L_n(p)|k_n))_3$$
for all $n \geq n_0$, i.e.
$$G(L_\infty(p)|k_n)/(G_\infty)_3 = G_\infty/(G_\infty)_3 \rtimes \Gamma_n \cong G(L_n(p)|k_n)/(G(L_n(p)|k_n))_3 \rtimes \Gamma_n.$$

Then $\langle \gamma_j^{p^{n-n_1}} \rangle$ acts trivially on $G_\infty/(G_\infty)_3$ for all $j \leq s$ and $N_n$ is contained in $(G_\infty)_3$.

Suppose that $G(L_n(p)|k_n)$, $n \geq n_0$, is powerful. Then
$$[G_\infty, G_\infty] \subseteq (G_\infty)^p N_n.$$

By assumption on $n_0$ the group $N_n$ is contained in $(G_\infty)_3$, and so

$$[G_\infty, G_\infty] \subseteq (G_\infty)^p [G_\infty, [G_\infty, G_\infty]].$$

From this inclusion it follows that

$$[G_\infty, G_\infty] \subseteq (G_\infty)^p,$$

thus $G_\infty$ is powerful.

Using proposition 2.1, we can assume that

$$d_{k_n}^+ = \dim_{\mathbb{F}_p}(Cl(k_n)/p)^+ \geq 1.$$

Let $K|k_n$ be an unramified Galois extension of degree $p$ such that $G(K|k_n) = G(K|k_n)^+$. Because of our definition of $n_1$ the field $K$ is not contained in $k_\infty$ and $G(L_\infty(p)|K_\infty)$ is a normal subgroup of $G(L_\infty(p)|k_\infty)$ of index $p$. Using results of Iwasawa theory, [6] (11.4.13) and (11.4.8), we obtain

$$d(G(L_\infty(p)|K_\infty))^- = p(d(G(L_\infty(p)|k_\infty))^- - 1) + 1.$$

From [2] theorem 3.8 and the equality above it follows that

$$\begin{aligned}
& d(G(L_\infty(p)|k_\infty))^+ + d(G(L_\infty(p)|k_\infty))^- \\
= & \; d(G(L_\infty(p)|k_\infty)) \\
\geq & \; d(G(L_\infty(p)|K_\infty)) \\
= & \; d(G(L_\infty(p)|K_\infty))^+ + d(G(L_\infty(p)|K_\infty))^- \\
= & \; d(G(L_\infty(p)|K_\infty))^+ + p(d(G(L_\infty(p)|k_\infty))^- - 1) + 1.
\end{aligned}$$

The maximal quotient $G(L_\infty(p)|k_\infty)_\Delta$ of $G(L_\infty(p)|k_\infty)$ with trivial action of $\Delta$ is also powerful and we have $d(G(L_\infty(p)|k_\infty)_\Delta) = d(G(L_\infty(p)|k_\infty))^+$. Using again [2] theorem 3.8, we get

$$d(G(L_\infty(p)|k_\infty))^+ \geq d(G(L_\infty(p)|K_\infty))^+.$$

Both inequalities together imply

$$d(G(L_\infty(p)|k_\infty))^- \leq 1.$$

Using [6] (11.4.4), we finally obtain

$$d(G(L_\infty(p)|k_\infty))^+, \; d(G(L_\infty(p)|k_\infty))^- \leq 1.$$

It follows that $G(L_n(p)|k_n)$ is a powerful pro-$p$-group with $d(G(L_n(p)|k_n)) \leq 2$. If $G(L_n(p)|k_n)$ is not finite, then it contains an open subgroup $P$ which is a Poincaré group (see [5] chap.V theorem (2.2.8) and (2.5.8)) of dimension $\dim P = d(P) \leq 2$ (use again [2] theorem 3.8). But these groups have the group $\mathbb{Z}_p$ as homomorphic image. By the finiteness of the class number it follows that $G(L_n(p)|k_n)$ is finite. □

**Remark:** The theorems 3.1, 3.2 and 3.3 above hold, if we replace $L(p)$ by $L_S(p)$ and $Cl$ by $Cl_S$ where $S \supseteq S_\infty$ is a set of primes which do not split in the extension $k|k^+$.

Now we consider the conjecture for general $p$-adic analytic groups. Let

$$1 \longrightarrow \mathcal{D} \longrightarrow \mathcal{G} \longrightarrow G \longrightarrow 1$$

be an exact sequence of pro-$p$-groups. For an open normal subgroup $H$ of $G$ we denote the pre-image of $H$ in $\mathcal{G}$ by $\mathcal{H}$. Thus we get a commutative exact diagram

$$
\begin{array}{ccccccccc}
1 & \longrightarrow & \mathcal{D} & \longrightarrow & \mathcal{G} & \longrightarrow & G & \longrightarrow & 1 \\
& & \| & & \cup\!\!\!| & & \cup\!\!\!| & & \\
1 & \longrightarrow & \mathcal{D} & \longrightarrow & \mathcal{H} & \longrightarrow & H & \longrightarrow & 1.
\end{array}
$$

**Proposition 3.4** *With the notation as above assume that*
  (i) *$\mathcal{G}$ is finitely generated and $cd_p \mathcal{G} \leq 2$,*
 (ii) *$cd_p G < \infty$,*
(iii) *the Euler-Poincaré characteristic of $\mathcal{G}$ is zero, i.e.*

$$\chi(\mathcal{G}) = \sum_{i=0}^{2}(-1)^i \dim_{\mathbb{F}_p} H^i(\mathcal{G}, \mathbb{Z}/p\mathbb{Z}) = 0 .$$

*Then*
  $d(\mathcal{H})$ *is unbounded for varying open normal subgroups $H$ of $G$ or $cd_p G \leq 2$.*

**Proof:** Suppose that $\dim_{\mathbb{F}_p} H^1(\mathcal{H}, \mathbb{Z}/p\mathbb{Z})$ is bounded for varying $H$. Since $\chi(\mathcal{G}) = 0$, the same is true for $\dim_{\mathbb{F}_p} H^2(\mathcal{H}, \mathbb{Z}/p\mathbb{Z})$. It follows that $H^i(\mathcal{D}, \mathbb{Z}/p\mathbb{Z})$ is finite for $i = 1, 2$. By [6] proposition (3.3.7), we obtain

$$cd_p \mathcal{G} = cd_p G + cd_p \mathcal{D} \geq cd_p G.$$

This proves the proposition. □

As an application to our problem we get the following result for the maximal unramified $p$-extension $L_k(p)$ of a number field $k$.

**Theorem 3.5** *Let $p \neq 2$ and let $k$ be a CM-field containing $\mu_p$ with maximal totally real subfield $k^+$. Assume that $\mu_p \not\subseteq k^+_{\mathfrak{p}}$ for all primes $\mathfrak{p}$ of $k^+$ above $p$. Then the following holds:*

$$\begin{array}{lll} either & (i) & G(L_{k^+}(p)|k^+) \text{ is finite,} \\ or & (ii) & G(L_k(p)|k) \quad \text{ is not } p\text{-adic analytic,} \end{array}$$

*with other words, if $G(L_k(p)|k)$ is p-adic analytic, then $G(L_{k^+}(p)|k^+)$ is finite.*

**Proof:** Suppose that (i) and (ii) do not hold. Then the maximal quotient $G(L_{k^+}(p)|k^+)$ of the $p$-adic analytic group $G(L_k(p)|k)$ with trivial action by $\Delta = G(k|k^+)$ is an infinite analytic group. Passing to a finite extension of $k^+$, we may assume that $G(L_{k^+}(p)|k^+)$ is uniform (our assumptions on $k$ are still valid). The dimension of $G(L_{k^+}(p)|k^+)$ is greater or equal to 3, since otherwise it would have the group $\mathbb{Z}_p$ as quotient which is impossible by the finiteness of the class number.

If $k^+_{S_p}(p)$ is the maximal $p$-extension of $k^+$ which is unramified outside $p$, then $cd_p \, G(k^+_{S_p}(p)|k^+) \leq 2$ and $\chi(G(k^+_{S_p}(p)|k^+)) = 0$, see [6] (8.3.17), (8.6.16) and (10.4.8). Applying proposition 3.4, we obtain that

$$\dim_{\mathbb{F}_p} H^1(G(k^+_{S_p}(p)|K^+), \mathbb{Z}/p\mathbb{Z}) = \dim_{\mathbb{F}_p} H^1(G(k_{S_p}(p)|K^+(\mu_p)), \mathbb{Z}/p\mathbb{Z})^+$$

is unbounded, if $K^+$ varies over the finite Galois extension of $k^+$ inside $L_{k^+}(p)$. By [6] theorem (8.7.3) and the assumption that $\mu_p \not\subseteq k^+_{\mathfrak{p}}$ for all primes $\mathfrak{p}|p$, it follows that

$$\begin{aligned} \dim_{\mathbb{F}_p} Cl(K^+(\mu_p))/p & \geq & \dim_{\mathbb{F}_p}(Cl_{S_p}(K^+(\mu_p))/p)^- \\ & = & \dim_{\mathbb{F}_p} H^1(G(k_{S_p}(p)|K^+(\mu_p)), \mathbb{Z}/p\mathbb{Z})^+ - 1 \end{aligned}$$

is unbounded for varying $K^+$ inside $L_{k^+}(p)$ and therefore $G(L_k(p)|k)$ is not $p$-adic analytic. This contradiction proves the theorem. $\square$

Mathematisches Institut
der Universität Heidelberg
Im Neuenheimer Feld 288
69120 Heidelberg
Germany

e-mail: wingberg@mathi.uni-heidelberg.de